# Lévy-based growth models

KRISTJANA ÝR JÓNSDÓTTIR, JÜRGEN SCHMIEGEL
and EVA B. VEDEL JENSEN

*T.N. Thiele Centre for Applied Mathematics in Natural Science, University of Aarhus, Denmark.
E-mail: eva@imf.au.dk*

In the present paper, we give a condensed review, for the nonspecialist reader, of a new modelling framework for spatio-temporal processes, based on Lévy theory. We show the potential of the approach in stochastic geometry and spatial statistics by studying Lévy-based growth modelling of planar objects. The growth models considered are spatio-temporal stochastic processes on the circle. As a by product, flexible new models for space–time covariance functions on the circle are provided. An application of the Lévy-based growth models to tumour growth is discussed.

## 1. Introduction

Stochastic spatio-temporal modelling is of great importance in a variety of disciplines of natural science, including biology (Cantalapiedra *et al.* (2001), Brix and Chadoeuf (2002), Fewster (2003), Gratzer *et al.* (2004)), image analysis (Feideropoulou and Pesquet-Popescu (2004)), geophysics (Calder (1986), Lovejoy *et al.* (1992), Sornette and Ouillon (2005)) and turbulence (Schmiegel *et al.* (2004), Schmiegel *et al.* (2005)), to name just a few. In particular, the modelling of tumour growth dynamics has been a very active research area in recent years (Delsanto *et al.* (2000), Peirolo and Scalerandi (2004), Pang and Tzeng (2004), Schmiegel (2006)). In most of the above-cited works, the model is given implicitly and the resulting dynamics are difficult to control explicitly. However, for applications and for the theoretical understanding of the modelling framework being employed it is essential to connect the ingredients of the model with dynamical and spatial properties of the system under consideration. Furthermore, for a parsimonious description of systems which are different with respect to the dynamics and physical mechanisms underlying the dynamics, it is desirable to have access to a flexible and, at the same time, mathematically tractable modelling framework.

Lévy-based models provide a promising modelling framework to meet these requirements concerning flexibility and dynamical control. Until now, Lévy-based models have mainly been used for describing turbulent flows (Barndorff-Nielsen and Schmiegel (2004),

Schmiegel *et al.* (2004), Schmiegel *et al.* (2005)). In the present paper, we show that Lévy-based spatio-temporal modelling also has important applications in stochastic geometry and spatial statistics. The main focus is on Lévy-based growth models, but we will also briefly touch upon another emerging area of application, that is, Lévy-driven Cox processes. We expect that the Lévy-based approach will have many applications in stochastic geometry and spatial statistics.

Lévy-based spatio-temporal models are constructed from Lévy bases, that is, infinitely divisible and independently scattered random measures. This terminology was recently introduced in Barndorff-Nielsen and Schmiegel (2004). One example of such a model describes the growth of a planar star-shaped object, using its radial function $R_t(\phi)$ at time $t$ and angle $\phi$. Here, $R_t(\phi)$ is the distance from a reference point $z$ to the boundary of the object at time $t$ in the direction $\phi \in [-\pi, \pi)$; see also Figure 2 below. The time derivative of the radial function is of the form

$$\frac{\partial}{\partial t}R_t(\phi) = \mu_t(\phi) + \int_{A_t(\phi)} f_t(\xi; \phi) Z(\mathrm{d}\xi), \qquad \phi \in [-\pi, \pi),$$

where $Z$ is a Lévy basis, $A_t(\phi) \subseteq [-\pi, \pi) \times (-\infty, t]$ is a subset of the past of time $t$, a so-called *ambit set* (cf. Barndorff-Nielsen and Schmiegel (2004)), $f_t(\cdot; \phi)$ is a deterministic weight function and $\mu_t$ a deterministic function. (In Latin, ambitus means either (1) the bounds or limits of a place or district, or (2) a sphere of action, expression or influence.) The induced model for the radial function is of the same type. An important advantage of these models is that explicit expressions for

$$\mathrm{Cov}(R_{t_1}(\phi_1), R_{t_2}(\phi_2))$$

can be derived in terms of the components of the model. This part of the paper is a natural continuation of the work initiated in Schmiegel (2006), which was mainly directed toward an audience of physicists. An introduction to different growth modelling approaches, including a short treatment of the Lévy-based approach adopted in the present paper, may be found in the SemStat contribution Jensen *et al.* (2006).

The paper is organized as follows. Section 2 provides some background on Lévy bases, which is the essential component of the modelling approach. In Section 3, Lévy-based spatio-temporal models are reviewed, while Lévy-based growth models are studied in Section 4. Section 5 contains explicit results for the covariance functions. In Section 6, an application of the Lévy-based growth models to tumour growth is discussed, while problems for future research are collected in Section 7.

## 2. Lévy bases

This section provides a brief overview of the general theory of Lévy bases, in particular, the theory of integration with respect to a Lévy basis. For a more detailed exposition, see Barndorff-Nielsen and Schmiegel (2004). As mentioned in the Introduction, a Lévy basis is an infinitely divisible and independently scattered random measure. Comprehensive accounts of the theory of independently scattered random measures may be found in Kallenberg (1989), Rajput and Rosinski (1989) and Kwapien and Woyczynski (1992).

Let $\mathcal{R}$ be a Borel subset of $\mathbb{R}^d$, let $\mathcal{B} = \mathcal{B}(\mathcal{R})$ be the Borel subsets of $\mathcal{R}$ and let $\mathcal{B}_b = \mathcal{B}_b(\mathcal{R})$ denote the class of bounded elements of $\mathcal{B}$. A collection of random variables $Z = \{Z(A) : A \in \mathcal{B}\}$ or $Z = \{Z(A) : A \in \mathcal{B}_b\}$ is said to be an *independently scattered random measure* if, for every sequence $\{A_n\}$ of disjoint sets in $\mathcal{B}$ (resp. $\mathcal{B}_b$), the random variables $Z(A_n)$, $n = 1, 2, \ldots$, are independent and $Z(\bigcup A_n) = \sum Z(A_n)$ a.s., where, in the case $Z = \{Z(A) : A \in \mathcal{B}_b\}$, we further require $\bigcup A_n \in \mathcal{B}_b$. We need to distinguish between the two cases $\mathcal{B}$ and $\mathcal{B}_b$ because the sums $\sum Z(A_n)$ must also be controlled in case where $Z$ can take both positive and negative values. If, moreover, $Z(A)$ is *infinitely divisible* for all $A \in \mathcal{B}$ or $\mathcal{B}_b$, $Z$ is called a *Lévy basis* (cf. Barndorff-Nielsen and Schmiegel (2004)).

For a random variable $X$, let us denote the cumulant function $\log \mathbb{E}(e^{i\lambda X})$ by $C(\lambda \ddagger X)$. When $Z$ is a Lévy basis, the cumulant function of $Z(A)$ can, by the famous Lévy–Khinchine representation, be written as

$$C(\lambda \ddagger Z(A)) = i\lambda a(A) - \tfrac{1}{2}\lambda^2 b(A) + \int_{\mathbb{R}} (e^{i\lambda u} - 1 - i\lambda u \mathbf{1}_{[-1,1]}(u)) U(du, A), \qquad (1)$$

where $a$ is a signed measure on $\mathcal{B}$ or $\mathcal{B}_b$, $b$ is a positive measure on $\mathcal{B}$ or $\mathcal{B}_b$ and $U(du, A)$ is a Lévy measure on $\mathbb{R}$ for fixed $A$ and a measure on $\mathcal{B}$ or $\mathcal{B}_b$ for fixed $du$. The measure $U$ will be referred to as the *generalized* Lévy measure. The Lévy basis $Z$ is said to have *characteristics* $(a, b, U)$. If $b = 0$, then $L$ is called a *Lévy jump basis* and if $U = 0$, then $L$ is a *Gaussian basis*; see the examples below. It follows from (1) that any Lévy basis $Z$ can be expressed as the sum of a Lévy jump basis $Z_1$ and an independent zero-mean Gaussian basis $Z_2$.

Without loss of generality (for details, see Rajput and Rosinski (1989)), we can assume that there exists a measure $\mu$ such that the generalized Lévy measure factorizes as

$$U(du, d\xi) = V(du, \xi)\mu(d\xi),$$

where $V(du, \xi)$ is a Lévy measure for fixed $\xi$. Furthermore, the measures $a$ and $b$ are absolutely continuous with respect to the measure $\mu$, that is,

$$a(d\xi) = \tilde{a}(\xi)\mu(d\xi), \qquad b(d\xi) = \tilde{b}(\xi)\mu(d\xi),$$

and $\tilde{a}$ and $\tilde{b}$ are uniformly bounded by some constant $C > 0$. One possible choice (but not the only one; see Rajput and Rosinski (1989)) for $\mu$ is

$$\mu(A) = |a|(A) + b(A) + \int_{\mathbb{R}} (1 \wedge r^2) U(dr, A),$$

where $|a|$ denotes the total variation measure of $a$ and $\wedge$ denotes minimum.

Let $Z'(\xi)$ be a random variable with the cumulant function

$$C(\lambda \ddagger Z'(\xi)) = i\lambda \tilde{a}(\xi) - \tfrac{1}{2}\lambda^2 \tilde{b}(\xi) + \int_{\mathbb{R}} (e^{i\lambda u} - 1 - i\lambda u \mathbf{1}_{[-1,1]}(u)) V(du, \xi).$$

Then,

$$C(\lambda \ddagger Z(\mathrm{d}\xi)) = C(\lambda \ddagger Z'(\xi))\mu(\mathrm{d}\xi). \qquad (2)$$

If $\tilde{a}(\xi)$, $\tilde{b}(\xi)$ and the Lévy measure $V(\cdot;\xi)$ do not depend on $\xi$, then we call $Z$ a factorizable Lévy basis and $Z'(\xi) = Z'$ also does not depend on $\xi$. If, moreover, $\mu$ is proportional to the Lebesgue measure, then $Z$ is called a *homogeneous* Lévy basis and all finite-dimensional distributions of $Z$ are translation invariant.

The usefulness of the above definitions becomes clear in connection with the integration of a measurable function $f$ on $\mathcal{R}$ with respect to a Lévy basis $Z$. For simplicity, we denote this integral by $f \cdot Z$. Important for many calculations is the following equation for the cumulant function of the stochastic integral $f \cdot Z$ (subject to minor regularity conditions, cf., for instance, Barndorff-Nielsen and Thorbjørnsen (2003))

$$C(\lambda \ddagger f \cdot Z) = \int C(\lambda f(\xi) \ddagger Z'(\xi))\mu(\mathrm{d}\xi). \qquad (3)$$

The result (3) can be heuristically derived from (2). A similar result holds for the logarithm of the Laplace transform of $f \cdot Z$ (assumed to be finite),

$$K(\lambda \ddagger f \cdot Z) = \int K(\lambda f(\xi) \ddagger Z'(\xi))\mu(\mathrm{d}\xi). \qquad (4)$$

The function $K$ will, in the following, be called the *kumulant function*.

We will now give a few examples of Lévy bases.

***Example 1 (Gaussian Lévy basis).*** If $Z$ is a Lévy basis with $Z(A) \sim N(a(A), b(A))$, where $a$ is a signed measure and $b$ is a positive measure, we call $Z$ a Gaussian Lévy basis. The *Gaussian* Lévy basis has characteristics $(a, b, 0)$ and its cumulant function is

$$C(\lambda \ddagger Z(A)) = \mathrm{i}\lambda a(A) - \tfrac{1}{2}\lambda^2 b(A).$$

We have $Z'(\xi) \sim N(\tilde{a}(\xi), \tilde{b}(\xi))$, that is, $C(\lambda \ddagger Z'(\xi)) = \mathrm{i}\lambda\tilde{a}(\xi) - \tfrac{1}{2}\lambda^2\tilde{b}(\xi)$. Furthermore,

$$C(\lambda \ddagger f \cdot Z) = \int C(\lambda f(\xi) \ddagger Z'(\xi))\mu(\mathrm{d}\xi) = \mathrm{i}\lambda(f \cdot a) - \tfrac{1}{2}\lambda^2(f^2 \cdot b). \qquad (5)$$

Note that $f \cdot Z \sim N(f \cdot a, f^2 \cdot b)$.

***Example 2 (Lévy jump basis).*** A Lévy basis is called a Lévy jump basis if the characteristics of the basis are $(a, 0, U)$. In Table 1, we specify the functions $V$ and $\tilde{a}$ for three important examples of Lévy jump bases: the Poisson basis, the Gamma basis and the inverse Gaussian basis. We also list the distribution of the random variable $Z'(\xi)$, its cumulant function, mean and variance. All parameters are positive.

Note that if $Z$ is a Poisson basis, then $Z(A) \sim \mathrm{Po}(\mu(A))$ with probability function

$$\frac{\mathrm{e}^{-\mu(A)}\mu(A)^x}{x!}, \qquad x = 0, 1, 2, \ldots.$$

**Table 1.** The definition of three Lévy jump bases – the Poisson basis, the Gamma basis and the inverse Gaussian basis – and the distribution of $Z'(\xi)$, with the corresponding cumulant function, mean and variance

|  | Poisson | Gamma | Inverse Gaussian |
|---|---|---|---|
| $V(\mathrm{d}u, \xi)$ | $\delta_1(\mathrm{d}u)$ | $\mathbf{1}_{\mathbb{R}_+}(u)\beta u^{-1}\mathrm{e}^{-\alpha(\xi)u}\,\mathrm{d}u$ | $\frac{\eta}{\sqrt{2\pi}}\mathbf{1}_{\mathbb{R}_+}(u)u^{-3/2}\mathrm{e}^{(-1/2)\gamma^2(\xi)u}\,\mathrm{d}u$ |
| $\tilde{a}(\xi)$ | 1 | $\beta\left(\frac{1-\mathrm{e}^{-\alpha(\xi)}}{\alpha(\xi)}\right)$ | $\frac{\eta}{\sqrt{2\pi}}\int_0^1 u^{-1/2}\mathrm{e}^{-(1/2)\gamma^2(\xi)u}\,\mathrm{d}u$ |
| $Z'(\xi)$ | Po(1) | $\Gamma(\beta, \alpha(\xi))$ | $\mathrm{IG}(\eta, \gamma(\xi))$ |
| $C(\lambda \ddagger Z'(\xi))$ | $\mathrm{e}^{\mathrm{i}\lambda} - 1$ | $-\beta\log\left(1 - \frac{\mathrm{i}\lambda}{\alpha(\xi)}\right)$ | $\eta\gamma(\xi)\left(1 - \sqrt{1 - \frac{2\mathrm{i}\lambda}{\gamma^2(\xi)}}\right)$ |
| $\mathbb{E}(Z'(\xi))$ | 1 | $\frac{\beta}{\alpha(\xi)}$ | $\frac{\eta}{\gamma(\xi)}$ |
| $\mathbb{V}(Z'(\xi))$ | 1 | $\frac{\beta}{\alpha^2(\xi)}$ | $\frac{\eta}{\gamma^3(\xi)}$ |

If $Z$ is a Gamma basis with $\alpha(\xi) \equiv \alpha$, then $Z(A) \sim \Gamma(\beta\mu(A), \alpha)$ with density

$$\frac{\alpha^{\beta\mu(A)}}{\Gamma(\beta\mu(A))}x^{\beta\mu(A)-1}\mathrm{e}^{-\alpha x}, \qquad x > 0,$$

while if $Z$ is an inverse Gaussian basis with $\gamma(\xi) \equiv \gamma$, then $Z(A) \sim \mathrm{IG}(\eta\mu(\mathrm{A}), \gamma)$ with density

$$\frac{\eta\mu(A)\mathrm{e}^{\eta\mu(A)\gamma}}{\sqrt{2\pi}}x^{-3/2}\exp\left\{-\frac{1}{2}((\eta\mu(A))^2 x^{-1} + \gamma^2 x)\right\}, \qquad x > 0.$$

The Poisson, Gamma and inverse Gaussian Lévy bases are examples of the random $G$-measures introduced in Brix (1999). These measures are purely discrete and can be written as (the Lévy–Itô representation)

$$Z(A) = a_0(A) + \int_{\mathbb{R}_+} xN(\mathrm{d}x, A), \qquad (6)$$

where $N$ is a Poisson basis on $\mathbb{R}_+ \times \mathcal{R}$ with intensity measure $U$ and

$$a_0(A) = a(A) - \int_0^1 xU(\mathrm{d}x, A).$$

Note that equation (6) can also be written as

$$Z(A) = a_0(A) + \sum_{(u,\xi)\in\Phi} u\mathbf{1}_A(\xi), \qquad (7)$$

where $\Phi$ is a Poisson point process on $\mathbb{R}_+ \times \mathcal{R}$ with intensity function $U$. If $f$ is a measurable function on $\mathcal{R}$, we have

$$f \cdot Z = f \cdot a_0 + \sum_{(u,\xi) \in \Phi} u f(\xi).$$

Finally, it should be noted that any Lévy process $\{Z_t\}_{t \in \mathbb{R}}$ induces a Lévy basis $Z$ on $\mathcal{R}$ by

$$Z((a,b]) = Z_b - Z_a, \qquad a, b \in \mathbb{R}.$$

## 3. Lévy-based spatio-temporal modelling

Let us consider a random variable $X_t(\sigma)$ depending on time $t$ and a position $\sigma$ in space. In the following, we will assume that $(\sigma, t) \in \mathcal{R} = \mathcal{S} \times \mathbb{R}$, where $\mathcal{S} \subseteq \mathbb{R}^n$, say. A Lévy-based spatio-temporal model for $X = \{X_t(\sigma) : (\sigma, t) \in \mathcal{R}\}$ is based on the intuitive picture of an ambit set $A_t(\sigma)$ associated with each point $(\sigma, t) \in \mathcal{R}$, which defines the dependency on the past at time $t$ and position $\sigma$. The ambit set $A_t(\sigma)$ will always satisfy the following conditions:

$$(\sigma, t) \in A_t(\sigma),$$
$$A_t(\sigma) \subseteq \mathcal{S} \times (-\infty, t].$$

An illustration is shown in Figure 1. The *linear spatio-temporal Lévy model* for $X = \{X_t(\sigma) : (\sigma, t) \in \mathcal{R}\}$ is then defined as

$$X_t(\sigma) = \int_{A_t(\sigma)} f_t(\xi; \sigma) Z(\mathrm{d}\xi), \qquad (8)$$

where $Z$ is a Lévy basis and $f_t(\xi; \sigma)$ is a deterministic weight function, which is assumed to be suitable for the integral to exist. The process

$$\tilde{X} = \{\exp(X_t(\sigma)) : (\sigma, t) \in \mathcal{R}\}$$

is said to follow an *exponential spatio-temporal Lévy model*.

The spatio-temporal Lévy models can be viewed as generalizations of the familiar moving average processes in time series, extended (i) from discrete to continuous time, (ii) from one dimension (time) to space-time, preserving a notion of causality in time (i.e., the ambit set lies in the past of time $t$) and (iii) from Gaussian white noise to more general infinitely divisible processes.

Using the key relation (3), we can derive expressions for moments in the linear spatio-temporal Lévy model. Thus, we find

$$\mathbb{E}(X_t(\sigma)) = \int_{A_t(\sigma)} f_t(\xi; \sigma) \mathbb{E}(Z'(\xi)) \mu(\mathrm{d}\xi) \qquad (9)$$

and

$$\mathbb{V}(X_t(\sigma)) = \int_{A_t(\sigma)} f_t^2(\xi;\sigma)\mathbb{V}(Z'(\xi))\mu(\mathrm{d}\xi), \tag{10}$$

where $\mathbb{V}$ is the notation used for variance. The covariances are of the form

$$\mathrm{Cov}(X_{t_1}(\sigma_1), X_{t_2}(\sigma_2)) = \int_{A_{t_1}(\sigma_1) \cap A_{t_2}(\sigma_2)} f_{t_1}(\xi;\sigma_1) f_{t_2}(\xi;\sigma_2)\mathbb{V}(Z'(\xi))\mu(\mathrm{d}\xi). \tag{11}$$

If the weight function is constant, that is $f_t(\xi;\sigma) \equiv f$, and if the Lévy basis $Z$ is factorizable, then (11) reduces to

$$\mathrm{Cov}(X_{t_1}(\sigma_1), X_{t_2}(\sigma_2)) = f^2\mathbb{V}(Z')\mu(A_{t_1}(\sigma_1) \cap A_{t_2}(\sigma_2)). \tag{12}$$

In this case, the covariance structure depends only on the $\mu$-measure of the intersection of the two ambit sets.

Equation (4) enables us to calculate arbitrary mixed $n$th order moments of $\tilde{X}_t(\sigma) = \exp(X_t(\sigma))$. If the moments are finite, then

$$\mathbb{E}(\tilde{X}_{t_1}(\sigma_1) \cdots \tilde{X}_{t_n}(\sigma_n)) = \exp\left(\int_{\mathcal{R}} K\left(\sum_{j=1}^n f_{t_j}(\xi;\sigma_j)\mathbf{1}_{A_{t_j}(\sigma_j)}(\xi) \ddagger Z'(\xi)\right)\mu(\mathrm{d}\xi)\right). \tag{13}$$

The corresponding expressions for the mixed $n$th order moments of $X_t(\sigma)$ are obtained from

$$\mathbb{E}(X_{t_1}(\sigma_1) \cdots X_{t_n}(\sigma_n)) = \frac{\partial^n}{\partial \lambda_1 \cdots \partial \lambda_n} \mathbb{E}(\tilde{X}_{t_1}^{\lambda_1}(\sigma_1) \cdots \tilde{X}_{t_n}^{\lambda_n}(\sigma_n))\bigg|_{\lambda_1=\cdots=\lambda_n=0}, \tag{14}$$

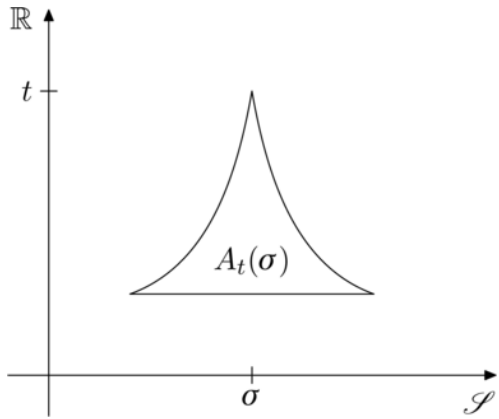

**Figure 1.** The ambit set $A_t(\sigma)$.

where

$$\mathbb{E}(\tilde{X}_{t_1}^{\lambda_1}(\sigma_1) \cdots \tilde{X}_{t_n}^{\lambda_n}(\sigma_n)) = \exp\left(\int_{\mathcal{R}} K\left(\sum_{j=1}^{n} \lambda_j f_{t_j}(\xi;\sigma_j)\mathbf{1}_{A_{t_j}(\sigma_j)}(\xi) \ddagger Z'(\xi)\right)\mu(\mathrm{d}\xi)\right). \tag{15}$$

The relative second-order moments of $\tilde{X}_t(\sigma)$ have a particularly attractive form

$$\frac{\mathbb{E}(\tilde{X}_{t_1}(\sigma_1)\tilde{X}_{t_2}(\sigma_2))}{\mathbb{E}(\tilde{X}_{t_1}(\sigma_1))\mathbb{E}(\tilde{X}_{t_2}(\sigma_2))} = \exp\left(\int_{A_{t_1}(\sigma_1) \cap A_{t_2}(\sigma_2)} g(\xi;t_1,t_2,\sigma_1,\sigma_2)\mu(\mathrm{d}\xi)\right), \tag{16}$$

where

$$g(\xi;t_1,t_2,\sigma_1,\sigma_2) \\ = K((f_{t_1}(\xi;\sigma_1) + f_{t_2}(\xi;\sigma_2)) \ddagger Z'(\xi)) - K(f_{t_1}(\xi;\sigma_1) \ddagger Z'(\xi)) - K(f_{t_2}(\xi;\sigma_2) \ddagger Z'(\xi)).$$

In the simple case where the weight functions are constant, that is, $f_t(\xi;\sigma) \equiv f$ for all $(\sigma,t) \in \mathcal{R}$ and $\xi \in \mathcal{R}$, and where the underlying Lévy basis is factorisable, $Z'(\xi) = Z'$, (16) reduces to

$$\exp(\overline{C}\mu(A_{t_1}(\sigma_1) \cap A_{t_2}(\sigma_2))), \tag{17}$$

where $\overline{C} = K(2f \ddagger Z') - 2K(f \ddagger Z')$. For a factorisable Lévy basis $Z$ and a constant weight function, one can express

$$\frac{\mathbb{E}(\tilde{X}_{t_1}^{\lambda_1}(\sigma_1) \cdots \tilde{X}_{t_n}^{\lambda_n}(\sigma_n))}{\mathbb{E}(\tilde{X}_{t_1}^{\lambda_1}(\sigma_1)) \cdots \mathbb{E}(\tilde{X}_{t_n}^{\lambda_n}(\sigma_n))} \tag{18}$$

in terms of different overlaps of the corresponding ambit sets (Schmiegel *et al.* (2005)).

## 4. Lévy-based growth models

In this section, we demonstrate the potential of the Lévy setup in stochastic geometry and spatial statistics by constructing Lévy-based stochastic models for growing objects. We focus on planar objects, but generalisations to higher dimensions are straightforward. We denote the planar object at time $t$ by $Y_t \subset \mathbb{R}^2$ and will assume that $Y_t$ is compact and star-shaped with respect to a point $z \in Y_t$ for all $t$. The boundary of the star-shaped object $Y_t$ can be determined by its radial function $R_t = \{R_t(\phi) : \phi \in [-\pi,\pi)\}$, where

$$R_t(\phi) = \max\{r : z + r(\cos\phi,\sin\phi) \in Y_t\}, \qquad \phi \in [-\pi,\pi)$$

(cf. Figure 2).

The growth rate will be described by the equation

$$\frac{\partial}{\partial t}R_t(\phi) = \mu_t(\phi) + \int_{A_t(\phi)} f_t(\xi;\phi)Z(\mathrm{d}\xi). \tag{19}$$

Here, the deterministic function $\mu_t : [-\pi, \pi) \to \mathbb{R}$ contributes to the overall growth pattern while the stochastic integral determines the dependence structure in the growth process. The ambit set $A_t(\phi) \subseteq [-\pi, \pi) \times (-\infty, t]$ relates to past events, $f_t(\cdot; \phi) : [-\pi, \pi) \times \mathbb{R} \to \mathbb{R}$ is a deterministic weight function (assumed to be suitable for the integral to exist) and $Z$ is a Lévy basis on $[-\pi, \pi) \times \mathbb{R}$. The weight functions and ambit sets must be defined cyclically in the angle such that the radial function $R_t(\phi)$ becomes cyclic. In the following, all angular calculations are regarded as cyclic.

Note that for nonnegative weight functions and nonnegative Lévy bases, $\frac{\partial}{\partial t} R_t(\phi) \geq 0$. In other cases, equation (19) still has a growth interpretation if the right-hand side of (19) is non-negative with probability one.

Using (9), the mean growth rate becomes

$$\mathbb{E}\left(\frac{\partial}{\partial t} R_t(\phi)\right) = \mu_t(\phi) + \int_{A_t(\phi)} f_t(\xi; \phi) \mathbb{E}(Z'(\xi)) \mu(\mathrm{d}\xi).$$

In the special case where $Z$ is a zero-mean Gaussian Lévy basis, $\mu_t(\phi)$ is indeed the mean growth rate at time $t$ in direction $\phi$. In other cases, $\mu_t(\phi)$ must be chosen such that the mean growth rate becomes as desired. There is a large literature on deterministic modelling of growth. A classical example is the Gompertz growth rate specified by

$$\mathbb{E}\left(\frac{\partial}{\partial t} R_t(\phi)\right) = \mu_t = \kappa_0 \exp\left[\frac{\eta}{\gamma}(1 - \exp(-\gamma t))\right] \eta \exp(-\gamma t)$$

(cf., e.g., Steel (1977)).

The ambit set $A_t(\phi)$ plays an important role in this modelling approach and affects the degree of dependence on the past. The extent of the dependence on the past may be specified by the minimal time-lag $T(t)$ such that

$$A_t(\phi) \subseteq [-\pi, \pi) \times [t - T(t), t], \qquad \phi \in [-\pi, \pi).$$

For an illustration, see Figure 3. Note that it follows from the fact that $Z$ is independently scattered that the random growth rates at time $t_1$ and $t_2$ are independent if $\min(t_1, t_2) < \max(t_1 - T(t_1), t_2 - T(t_2))$. The actual form of the ambit set $A_t(\phi)$ will depend on the specific growth process being modelled. A number of examples are given below. The induced correlation structure will be discussed in more detail in Section 5. A discrete version of (19) with a Gaussian Lévy basis has earlier been discussed in Jónsdóttir and Jensen (2005).

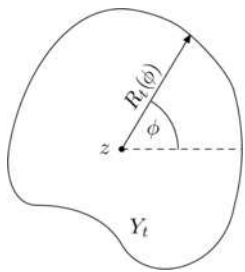

**Figure 2.** The star-shaped object $Y_t$ is determined by its radial function $R_t(\phi)$ at time $t$ and angle $\phi$.

For the interpretation of (19) as a growth model, it is helpful to represent the ambit set as a stochastic subset of the growing object. This is possible if the stochastic time transformation $t \to R_t(\phi)$ is non-decreasing for each $\phi \in [-\pi, \pi)$. We can then represent the ambit set $A_t(\phi)$ as a subset of $Y_t$,

$$\tilde{A}_t(\phi) = \{(R_s(\theta)\cos\theta, R_s(\theta)\sin\theta) : (\theta, s) \in A_t(\phi)\}.$$

It follows from the fact that $A_t(\phi) \subseteq [-\pi, \pi) \times (-\infty, t]$ that $\tilde{A}_t(\phi)$ is actually a subset of $Y_t$. Furthermore, since $(\phi, t) \in A_t(\phi)$, the set $\tilde{A}_t(\phi)$ touches the boundary of $Y_t$ at the point $(R_t(\phi)\cos\phi, R_t(\phi)\sin\phi)$. It is the 'events' in $\tilde{A}_t(\phi)$ that influence the growth rate at time $t$ in direction $\phi$. Figure 4 illustrates the set $\tilde{A}_t(\phi)$.

In the particular case where $Z$ is a Poisson basis and $\Psi$ the associated Poisson point process on $[-\pi, \pi) \times \mathbb{R}$, we can represent that part of the spatio-temporal point process $\Psi$, arrived before time $t$,

$$\Psi_t = \{(\theta_i, t_i) : t_i \leq t\},$$

as a subset of $Y_t$:

$$\tilde{\Psi}_t = \{(R_{t_i}(\theta_i)\cos\theta_i, R_{t_i}(\theta_i)\sin\theta_i) : t_i \leq t\}.$$

We can think of $\tilde{\Psi}_t$ as consisting of locations of outbursts at time points before $t$. Finally, if we let

$$\tilde{f}_t((s\cos\theta, s\sin\theta); \phi) = f_t((\theta, s); \phi),$$

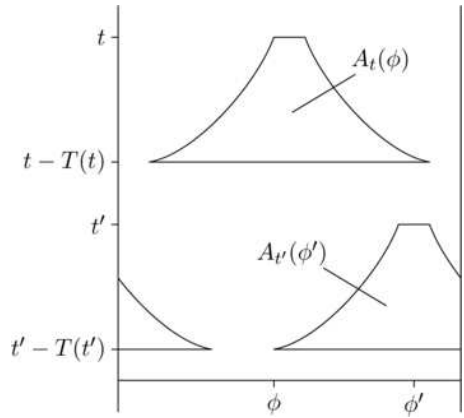

**Figure 3.** Two ambit sets $A_t(\phi)$ and $A_{t'}(\phi')$. Note the cyclic definition in the angle. The vertical lines are $\phi = -\pi$ and $\phi = \pi$, respectively.

then the fundamental equation (19) can be written as

$$\frac{\partial}{\partial t}R_t(\phi) = \mu_t(\phi) + \sum_{\tilde{\xi} \in \tilde{\Psi}_t \cap \tilde{A}_t(\phi)} \tilde{f}_t(\tilde{\xi}; \phi). \qquad (20)$$

According to (20), the growth rate at time $t$ in the direction $\phi$ depends on the outbursts at time points before $t$ which lie in the stochastic neighbourhood $\tilde{A}_t(\phi)$. This Poisson model is closely related to other recently suggested growth models (cf. Section 7.1).

Under (19), the induced model for $R_t(\phi)$ will be of the same linear form since

$$\begin{aligned} R_t(\phi) &= R_0(\phi) + \bar{\mu}_t(\phi) + \int_0^t \int_{A_s(\phi)} f_s(\xi; \phi) Z(\mathrm{d}\xi) \, \mathrm{d}s \\ &= R_0(\phi) + \bar{\mu}_t(\phi) + \int_{\bar{A}_t(\phi)} \bar{f}_t(\xi; \phi) Z(\mathrm{d}\xi), \end{aligned} \qquad (21)$$

where $R_0$ is the radial function at time $t=0$,

$$\bar{\mu}_t(\phi) = \int_0^t \mu_s(\phi) \, \mathrm{d}s,$$

$$\bar{A}_t(\phi) = \bigcup_{0 \le s \le t} A_s(\phi)$$

and

$$\bar{f}_t(\xi; \phi) = \int_0^t \mathbf{1}_{A_s(\phi)}(\xi) f_s(\xi; \phi) \, \mathrm{d}s. \qquad (22)$$

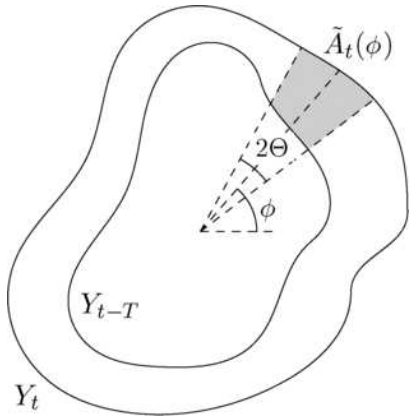

**Figure 4.** Illustration of the stochastic representation $\tilde{A}_t(\phi)$ (shown hatched) of the ambit set $A_t(\phi) = \{(\theta, s) : |\theta - \phi| \le \Theta, t - T \le s \le t\}$.

Note that the ambit sets associated with the radial function itself are increasing, that is,

$$t_1 \leq t_2 \Rightarrow \bar{A}_{t_1}(\phi) \subseteq \bar{A}_{t_2}(\phi).$$

If $t - T(t)$ is a non-decreasing function of $t$, then, because $Z$ is independently scattered,

$$R_{\leq t_1 - T(t_1)} = \{R_{t_2}(\phi) : t_2 \leq t_1 - T(t_1), \phi \in [-\pi, \pi)\}$$

and

$$R_{>t_1} - R_{t_1} = \{R_{t_2}(\phi) - R_{t_1}(\phi) : t_2 > t_1, \phi \in [-\pi, \pi)\},$$

are independent.

The representation (21) is, of course, not unique. If, in particular,

$$A_t(\phi) = B_t \cap C_\phi, \tag{23}$$

then

$$\bar{A}_t(\phi) = \bar{B}_t \cap C_\phi,$$

where

$$\bar{B}_t = \bigcup_{0 \leq s \leq t} B_s$$

and we may choose, instead of (22),

$$\bar{f}_t(\xi; \phi) = \int_0^t \mathbf{1}_{B_s}(\xi) f_s(\xi; \phi) \, \mathrm{d}s.$$

Note that the only difference between the two versions of $\bar{f}_t(\xi; \phi)$ is the indicator function $\mathbf{1}_{C_\phi}$.

In some cases, it might be more natural to formulate the model in terms of the time derivative of $\ln(R_t(\phi))$,

$$\frac{\partial}{\partial t}(\ln(R_t(\phi))) = \mu_t(\phi) + \int_{A_t(\phi)} f_t(\xi; \phi) Z(\mathrm{d}\xi).$$

In this case, the induced model is an exponential spatio-temporal Lévy model,

$$R_t(\phi) = R_0(\phi) \exp\left(\bar{\mu}_t(\phi) + \int_{\bar{A}_t(\phi)} \bar{f}_t(\xi; \phi) Z(\mathrm{d}\xi)\right).$$

Mixed moments of $R_t(\phi)$ can be derived using the results in Section 3.

The choices of Lévy basis $Z$, ambit sets $A_t(\phi)$, weight functions $f_t(\xi; \phi)$ and $\mu_t(\phi)$ completely determine the growth dynamics. These four ingredients can be chosen arbitrarily and independently, which results in a great variety of different growth dynamics. We will now give a number of examples.

***Example 3.*** Consider a Lévy growth model for the time derivative of the radial function

$$\frac{\partial}{\partial t} R_t(\phi) = Z(A_t(\phi)), \tag{24}$$

where $Z$ is a Poisson Lévy basis with intensity measure concentrated on $[-\pi, \pi) \times \mathbb{R}_+$ of the form

$$\mu(\mathrm{d}\xi) = g(s)\,\mathrm{d}s\,\mathrm{d}\theta, \qquad \xi = (\theta, s).$$

Note that the corresponding point process in the Euclidean plane

$$\{(s\cos\theta, s\sin\theta) : (\theta, s) \text{ is a support point of } Z\}$$

constitutes a Poisson point process with intensity measure

$$\tilde{\mu}(\mathrm{d}x) = \frac{g(\|x\|)}{\|x\|}\,\mathrm{d}x, \qquad x \in \mathbb{R}^2.$$

In particular, if $g(s) = as$, $a > 0$, then the Poisson point process in the plane is homogeneous.

Let the ambit sets be given by

$$A_t(\phi) = \left\{(\theta, s) : |\theta - \phi| \leq \frac{\Theta}{s}, \max(0, t - T) \leq s \leq t\right\}.$$

Represented as subsets of the Euclidean plane, they will, as $t \to \infty$, approach rectangles of side lengths $2\Theta$ and $T$. Note that we can write the ambit set as

$$A_t(\phi) = B_t \cap C_\phi,$$

where

$$B_t = \{(\theta, s) : \max(0, t - T) \leq s \leq t\}$$

and

$$C_\phi = \left\{(\theta, s) : s \geq \frac{\Theta}{\pi}, |\theta - \phi| \leq \frac{\Theta}{s}\right\} \cup \left\{(\theta, s) : 0 \leq s \leq \frac{\Theta}{\pi}\right\}.$$

The mean growth rate at time $t$ and in the direction $\phi$ is, for $t > T + \frac{\Theta}{\pi}$,

$$\mu(A_t(\phi)) = 2\Theta \int_{t-T}^{t} \frac{g(s)}{s}\,\mathrm{d}s.$$

If $g(s) = as$, $a > 0$, then the mean growth rate is constant. Figure 5 shows simulations of this model with constant mean growth rate.

***Example 4.*** The size of the ambit sets plays an important role in the control of the local and global fluctuations of the boundary of the object $Y_t$. As an example, let us consider a Lévy growth model of the form

$$R_t(\phi) = \mu_t + Z(A_t(\phi)), \qquad (25)$$

where

$$A_t(\phi) = \{(\theta, s) : |\theta - \phi| \leq \Theta(s), t - T(t) \leq s \leq t\}.$$

In Figure 6, simulations are shown under this model, using a normal Lévy basis with

$$Z(A) \sim N(0, \sigma^2 \mu(A))$$

and $\mu$ equal to the Lebesgue measure on $\mathcal{R}$. Note that $\mu(A_t(\phi))$ does not depend on $\phi$. The simulations are based on a discretization of $Z$ on a grid with $\Delta t = 1$ and $\Delta \phi = \frac{2\pi}{1000}$. The upper and lower row of Figure 6 show simulations for two choices of angular extension of the ambit set at three different time points. The angular extension of the ambit set is $\Theta(s) = \frac{\pi}{100}$ for the upper row and $\Theta(s) = \frac{\pi}{5}$ for the lower row. For the smaller angular extension, we observe localized fluctuations of the profiles, but the global appearance is circular. For the larger angular extension, the fluctuations are on a much larger scale and the global appearance is more variable.

***Example 5.*** In this example, we study a model similar to the one described in Example 4, but now with a Gamma Lévy basis. The model equation is

$$R_t(\phi) = \tilde{\mu}_t + Z(A_t(\phi)), \qquad (26)$$

where $A_t(\phi)$ is defined as in Example 4,

$$Z(A) \sim \Gamma(\beta \mu(A), \alpha),$$

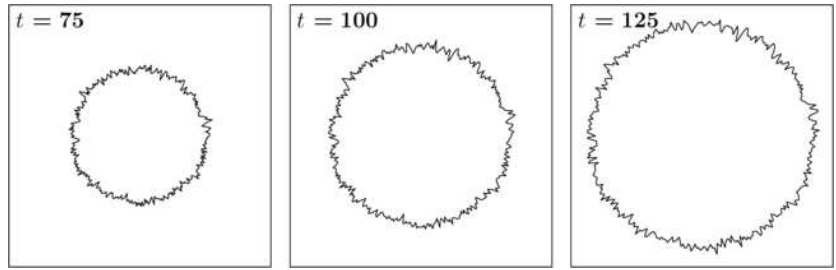

**Figure 5.** Simulation of the Lévy growth model (24) for the derivative of the radial function at time points $t = 75$, 100 and 125, using a Poisson Lévy basis. The parameters of the simulation are $g(s) = 10s$, $T = 1$ and $\Theta = 1/2$.

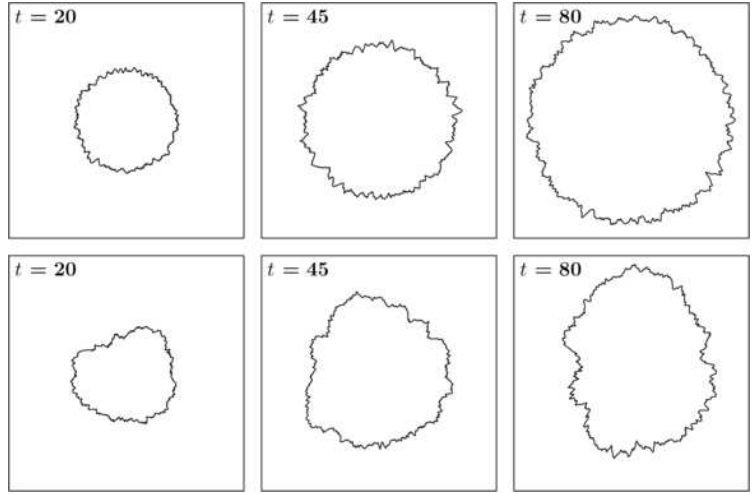

**Figure 6.** Simulation of the Lévy growth model (25) at time points $t = 20$, 45 and 80, using a Gaussian Lévy basis. The upper row and lower row show simulations of two choices of the angular extension of the ambit set $\Theta(s) = \frac{\pi}{100}$ and $\Theta(s) = \frac{\pi}{5}$, respectively. Otherwise, the parameters of the simulation are $\mu_{20} = 16$, $\mu_{45} = 24$, $\mu_{80} = 32$, $\sigma^2 = 1$ and $T(t) = t/5$.

and $\mu$ denotes Lebesgue measure on $\mathcal{R}$. The parameters $\alpha$, $\beta$ and $\tilde{\mu}_t$ are chosen such that $\mathbb{E}(R_t(\phi))$ and $\mathbb{V}(R_t(\phi))$ are the same as in the previous example. Here, we have used Table 1, together with (9) and (10). Accordingly, the parameters are chosen such that

$$\tilde{\mu}_t = \mu_t - \sigma\sqrt{\beta}\mu(A_t(0)),$$
$$\alpha = \sqrt{\frac{\beta}{\sigma^2}}.$$

The only free parameter is $\beta > 0$, which determines the skewness of the Gamma distribution of $Z(A_t(\phi))$. For large values of $\beta$, the distribution will resemble the Gaussian distribution.

The resulting simulations for $\beta = 1$ are shown in the upper and lower rows of Figure 7 for two choices of angular extension of the ambit set, $\Theta(s) = \frac{\pi}{100}$ and $\Theta(s) = \frac{\pi}{5}$, respectively. Note that more sudden outbursts are seen compared to the previous example.

**Example 6.** In Figure 8, we show simulations from the Lévy growth model

$$R_t(\phi) = f(\phi)(\mu_t + Z(A_t(\phi))), \tag{27}$$

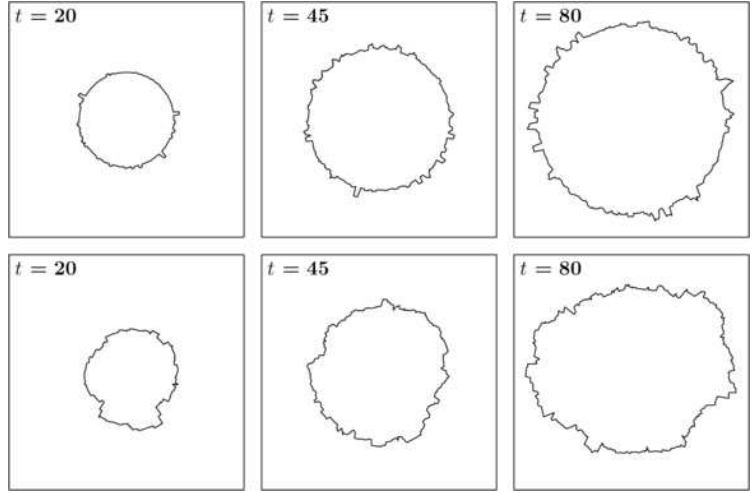

**Figure 7.** Simulation of the Lévy growth model (26) at time points $t = 20$, 45 and 80, using a Gamma Lévy basis. The upper row and lower row show simulations of two choices of the angular extension of the ambit set $\Theta(s) = \frac{\pi}{100}$ and $\Theta(s) = \frac{\pi}{5}$, respectively. Otherwise, $\beta = 1$ and the remaining parameters are determined by the parameters used in Example 4.

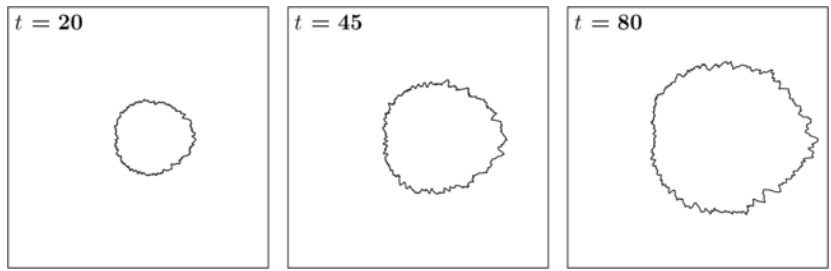

**Figure 8.** Simulation of the model (27) at time points $t = 20$, 45 and 80, using a Gaussian Lévy basis, with parameters as specified in Example 4. The weight function is given by $f_t(\phi) = 0.35 \exp(\frac{|\phi - \pi|}{\pi})$.

where $\mu_t$, $A_t(\phi)$ and $Z$ are specified as in Example 4 and

$$f_t(\phi) = 0.35 \exp\left(\frac{|\phi - \pi|}{\pi}\right).$$

Clearly, the growth of the object is asymmetric. The weight function $f_t(\phi)$ puts more weight on the angle $\phi_0 = 0$.

## 5. The induced covariance structure

Lévy-based growth models lead to flexible new models for space-time covariance functions on the circle, as we shall see in this section. The results presented here are of general interest for spatio-temporal processes on the circle.

We will derive expressions for $\mathrm{Cov}(R_{t_1}(\phi_1), R_{t_2}(\phi_2))$ under various assumptions on the Lévy basis $Z$, the ambit sets $A_t(\phi)$ and the weight functions $f_t(\xi; \phi)$. We will concentrate on the Lévy growth model (21) of linear form for $R_t$. Since we now are interested in covariances, it suffices to look at the model equation

$$R_t(\phi) = \int_{A_t(\phi)} f_t(\xi; \phi) Z(\mathrm{d}\xi),$$

where we, for simplicity, have omitted the bar on the ambit set and weight function. The covariance structure of $R_t(\phi)$ is then given by (cf. (11))

$$\mathrm{Cov}(R_{t_1}(\phi_1), R_{t_2}(\phi_2)) = \int_{A_{t_1}(\phi_1) \cap A_{t_2}(\phi_2)} f_{t_1}(\xi; \phi_1) f_{t_2}(\xi; \phi_2) \mathbb{V}(Z'(\xi)) \mu(\mathrm{d}\xi). \qquad (28)$$

Throughout this section, we will assume that

$$\begin{aligned} A_t(\phi) &= (\phi, 0) + A_t(0), \\ f_t(\xi; \phi) &= f_t((|\theta - \phi|, s); 0), \\ \mathbb{V}(Z'(\xi))\mu(\mathrm{d}\xi) &= g(s)\,\mathrm{d}s\,\mathrm{d}\theta \end{aligned} \qquad (29)$$

for all $\xi = (\theta, s) \in \mathcal{R}$ and $(\phi, t) \in \mathcal{R}$. These conditions ensure that $\mathrm{Cov}(R_{t_1}(\phi_1), R_{t_2}(\phi_2))$ only depends on the cyclic difference between $\phi_1$ and $\phi_2$. Accordingly, the spatio-temporal process

$$\{R_t(\phi) : t \in \mathbb{R}, \phi \in [-\pi, \pi)\}$$

will be second-order stationary in the space coordinate, but not necessarily in the time coordinate.

We will first consider the case where the angular extension of the ambit set is the full angular space, but the weight functions are arbitrary. Second, we consider the case of constant weight functions, but quite arbitrary ambit sets.

### 5.1. Ambit sets with full angular range

In this subsection, we consider ambit sets of the form

$$A_t(\phi) = [-\pi, \pi) \times [t - T(t), t].$$

In order to express the formulae as compactly as possible, we use, in the proposition below, the notation $t_1 \cap t_2$ for the time points shared by $A_{t_1}(\cdot)$ and $A_{t_2}(\cdot)$, that is,

$$t_1 \cap t_2 = \begin{cases} [\tilde{t}_1, \tilde{t}_2], & \text{if } \tilde{t}_1 \leq \tilde{t}_2, \\ \varnothing & \text{otherwise,} \end{cases}$$

where

$$\tilde{t}_1 = \max(t_1 - T(t_1), t_2 - T(t_2)) \quad \text{and} \quad \tilde{t}_2 = \min(t_1, t_2).$$

Using this notation, we can derive the following convenient and general expression for the covariances.

**Proposition 7.** *Let us assume that the ambit set is of the form $A_t(\phi) = [-\pi, \pi) \times [t - T(t), t]$ for all $(\phi, t) \in \mathcal{R}$ and let*

$$f_t(\xi; \phi) = a_0^t(s) + \sum_{k=1}^{\infty} a_k^t(s) \cos(k(\theta - \phi)), \tag{30}$$

$\xi = (\theta, s)$, *be the Fourier expansion of the weight function. The spatio-temporal covariances are then*

$$\text{Cov}(R_{t_1}(\phi_1), R_{t_2}(\phi_2)) = 2\tau_0(t_1, t_2) + \sum_{k=1}^{\infty} \tau_k(t_1, t_2) \cos(k(\phi_1 - \phi_2)), \tag{31}$$

*where*

$$\tau_k(t_1, t_2) = \pi \int_{t_1 \cap t_2} a_k^{t_1}(s) a_k^{t_2}(s) g(s) \, \mathrm{d}s.$$

**Proof.** The proof is straightforward. First, note that the actual form (30) of the Fourier expansion of the weight function is a consequence of (29). We obtain that

$$\text{Cov}(R_{t_1}(\phi_1), R_{t_2}(\phi_2))$$

$$= \int_{A_{t_1}(\phi_1) \cap A_{t_2}(\phi_2)} f_{t_1}(\xi; \phi_1) f_{t_2}(\xi; \phi_2) \mathbb{V}(Z'(\xi)) \mu(\mathrm{d}\xi)$$

$$= \pi \left[ 2 \int_{t_1 \cap t_2} a_0^{t_1}(s) a_0^{t_2}(s) g(s) \, \mathrm{d}s + \sum_{k=1}^{\infty} \left( \int_{t_1 \cap t_2} a_k^{t_1}(s) a_k^{t_2}(s) g(s) \, \mathrm{d}s \right) \cos(k(\phi_1 - \phi_2)) \right].$$

□

***Example 8.*** Suppose that the weight function is of the form (30) with $a_k^t(s) \equiv 0$ if $k \neq 1$. Then,

$$\text{Cov}(R_{t_1}(\phi_1), R_{t_2}(\phi_2)) = \pi \cos(\phi_1 - \phi_2) \int_{t_1 \cap t_2} a_1^{t_1}(s) a_1^{t_2}(s) g(s) \, \mathrm{d}s.$$

Since the covariance is a product of a spatial term and a temporal term, this model is separable (cf. Stein (2005) and references therein). The sign of the covariance may be positive or negative.

Note that, according to (31), the covariance $\text{Cov}(R_{t_1}(\phi_1), R_{t_2}(\phi_2))$ depends on $\phi_1$ and $\phi_2$ only via $|\phi_1 - \phi_2|$. For some choices of model parameters, the covariance also becomes stationary in the time coordinate. For instance, if $g(s) = 1$, $T(t) = T$ and $a_k^t(s) = b_k(t-s)$, then we have

$$\tau_k(t_1, t_2) = \pi \int_{\max(t_1-t_2,0)}^{T+\min(t_1-t_2,0)} b_k(u) b_k(t_2 - t_1 + u) \, du.$$

The induced model (31) for the covariance function is not, in general, separable in the sense that the covariance function can be written as a product of a term depending only on $t_1$ and $t_2$ and a term depending only on $\phi_1$ and $\phi_2$. This may be regarded as a strength of the model because separable covariance functions are often believed to give too simplistic a description of spatio-temporal data (cf., e.g., Stein (2005)). If, nevertheless, such simplifying assumptions are made, we obtain the following results.

**Corollary 9.** *Let the assumptions be as in Proposition 7. Assume that $a_k^t(s) = a_k^t$. The spatial correlations are then determined by the Fourier coefficients of the weight function $f$:*

$$\rho(R_t(\phi_1), R_t(\phi_2)) := \frac{\text{Cov}(R_t(\phi_1), R_t(\phi_2))}{\sqrt{\mathbb{V}(R_t(\phi_1)) \mathbb{V}(R_t(\phi_2))}} = \frac{2(a_0^t)^2 + \sum_{k=1}^{\infty} (a_k^t)^2 \cos(k(\phi_1 - \phi_2))}{2(a_0^t)^2 + \sum_{k=1}^{\infty} (a_k^t)^2}.$$

*If, in addition, $a_k^t = b_t c_k$, then the covariance model (31) is separable. Furthermore, the spatial correlations $\rho(R_t(\phi_1), R_t(\phi_2))$ do not depend on $t$, while the temporal correlations are determined by $T(t)$ and the function $g$:*

$$\rho(R_{t_1}(\phi), R_{t_2}(\phi)) = \frac{\int_{t_1 \cap t_2} g(s) \, ds}{[\int_{t_1-T(t_1)}^{t_1} g(s) \, ds \cdot \int_{t_2-T(t_2)}^{t_2} g(s) \, ds]^{1/2}}.$$

The covariance model (31) provides a possibility for extending stationary covariance functions on the circle (spatial covariance functions) to a spatio-temporal context. When $R_t$ is a stationary process on the circle, its covariance function can be expressed as

$$\text{Cov}(R_t(\phi_1), R_t(\phi_2)) = 2\lambda_0^t + \sum_{k=1}^{\infty} \lambda_k^t \cos(k(\phi_1 - \phi_2)). \tag{32}$$

Such a covariance function can be obtained by choosing the Fourier coefficients of the weight function to be

$$a_k^t(s) = a_k^t = \frac{1}{\sqrt{\pi}} \left[ \lambda_k^t \Big/ \int_{t-T(t)}^{t} g(s) \, ds \right]^{1/2}.$$

Note that there is still freedom in the modelling by choosing an arbitrary time-lag $T(t)$ and function $g$.

**Example 10.** The $p$th-order model for a stationary covariance function on the circle, described in Hobolth *et al.* (2003), has

$$\lambda_0^t = \lambda_1^t = 0, \lambda_k^t = [\alpha_t + \beta_t(k^{2p} - 2^{2p})]^{-1}, \qquad k = 2, 3, \ldots.$$

The model is called a $p$th-order model because it can be derived as a limit of discrete $p$th-order Markov models defined on a finite, systematic set of angles (cf. Hobolth and Jensen (2000)). This covariance structure is obtained by choosing

$$a_0^t(s) = a_1^t(s) = 0,$$

$$a_k^t(s) = \left[\pi \int_{t-T(t)}^{t} g(s)\, \mathrm{d}s\right]^{-1/2} [\alpha_t + \beta_t(k^{2p} - 2^{2p})]^{-1/2}, \qquad k = 2, 3, \ldots.$$

If $\alpha_t$ and $\beta_t$ are proportional, then the simplifying assumptions of Corollary 9 are fulfilled. In Jónsdóttir and Jensen (2005), this model has been used for the time derivative of the radial function. Only Gaussian Lévy bases are considered and neighbouring time points are assumed to be so far apart that the increments can be regarded as independent. The more general approach of the present paper allows for temporal correlations. Under the assumption $a_k^t = b_t c_k$, the temporal correlations are particularly simple. For instance, suppose that $T(t) \equiv 1$ and $t_2 - 1 \leq t_1 \leq t_2$. We then get for $g(s) = ae^{-bs}$, $a, b > 0$,

$$\rho(R_{t_1}(\phi), R_{t_2}(\phi)) = \frac{1}{\mathrm{e}^b - 1}[\mathrm{e}^{(1/2)b(t_1-t_2)+b} - \mathrm{e}^{(-1/2)b(t_1-t_2)}],$$

while, for $g(s) = as^\alpha$, $a > 0$, $\alpha \geq 1$,

$$\rho(R_{t_1}(\phi), R_{t_2}(\phi)) = \frac{t_1^{\alpha+1} - (t_2-1)^{\alpha+1}}{[(t_1^{\alpha+1} - (t_1-1)^{\alpha+1})(t_2^{\alpha+1} - (t_2-1)^{\alpha+1})]^{1/2}}.$$

Only in the first case are the temporal correlations always stationary.

## 5.2. Constant weight functions

In this subsection, we consider the case of constant weight functions. Without loss of generality, we assume that $f_t(\xi; \phi) \equiv 1$ and so (28) reduces to

$$\mathrm{Cov}(R_{t_1}(\phi_1), R_{t_2}(\phi_2)) = \int_{A_{t_1}(\phi_1) \cap A_{t_2}(\phi_2)} \mathbb{V}(Z'(\xi))\mu(\mathrm{d}\xi) = \mathbb{V}(Z')\mu(A_{t_1}(\phi_1) \cap A_{t_2}(\phi_2)), \tag{33}$$

where the last equality holds if the Lévy basis is factorisable.

It is not difficult (but sometimes tedious) to find explicit expressions for $\mathrm{Cov}(R_{t_1}(\phi_1), R_{t_2}(\phi_2))$ for specific choices of ambit sets. One simplifying assumption is to focus on ambit sets of the form

$$A_t(\phi) = B_t \cap C_\phi,$$

where

$$B_t = \{(\theta, s) : \max(0, t - T(t)) \leq s \leq t\},$$
$$C_\phi = \{(\theta, s) : |\phi - \theta| \leq \Theta(s)\}.$$

Usually, it is easier to find expressions for the temporal covariances than for the spatial covariances.

Evidently, (33) implies that $\mathrm{Cov}(R_{t_1}(\phi_1), R_{t_2}(\phi_2)) \geq 0$, which may be a severe restriction for the spatial covariances. In the proposition below, the spatial covariances are expressed in terms of the function delimiting the ambit set. The proposition provides insight into the class of spatial covariances that can be modelled using this approach.

**Proposition 11.** *Let $\mu(\mathrm{d}\xi) = g(s)\,\mathrm{d}s\,\mathrm{d}\theta$ for $\xi = (\theta, s)$. Let us suppose that there exists a continuous function $h_t : [-\pi, \pi) \to \mathbb{R}$ with the properties*

$$\begin{aligned} h_t(\theta) &= h_t(-\theta), \\ h_t \text{ is decreasing on } [0, \pi), \\ h_t(0) &= t \end{aligned} \tag{34}$$

*such that*

$$A_t(0) = \{(\theta, s) : h_t(\pi) \leq s \leq h_t(\theta)\}$$

*(cf. Figure 9). Let*

$$\bar{h}_t(\phi) = \int_0^{h_t(\phi)} g(s)\,\mathrm{d}s.$$

*If the Fourier expansion of $\bar{h}_t$ is ($\bar{h}_t(\phi) = \bar{h}_t(-\phi)$)*

$$\bar{h}_t(\phi) = \sum_{k=0}^{\infty} \gamma_k^t \cos(k\phi), \tag{35}$$

*then*

$$\mu(A_t(0) \cap A_t(\phi)) = \sum_{k=0}^{\infty} \lambda_k^t \cos(k\phi), \tag{36}$$

*where*

$$\lambda_0^t = \sum_{k\,\mathrm{odd}} \left[2\pi - \frac{16}{\pi k^2}\right] \gamma_k^t - 2\pi \sum_{k\,\mathrm{even}} \gamma_k^t,$$

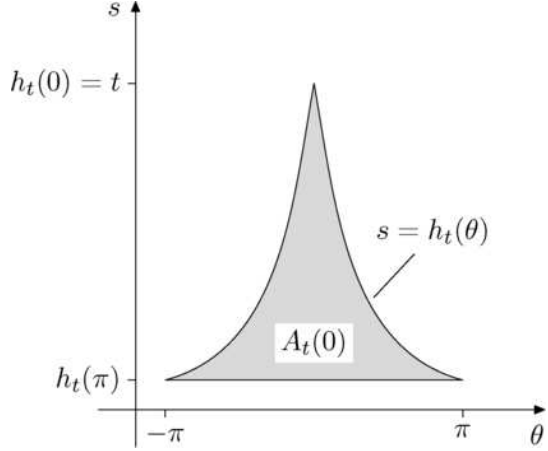

**Figure 9.** Illustration of the ambit set $A_t(0)$ bounded by the function $h_t$ (cf. (34)).

$$\lambda_j^t = \frac{16}{\pi} \sum_{k \text{ odd}} \frac{1}{(2j)^2 - k^2} \gamma_k^t, \qquad j = 1, 2, \ldots.$$

**Proof.** It is not difficult to show that

$$\mu(A_t(0) \cap A_t(\phi)) = 2 \int_{-\pi}^{-\pi+\phi/2} \bar{h}_t(\theta) \, d\theta + 2 \int_{\phi/2}^{\pi} \bar{h}_t(\theta) \, d\theta - 2\pi \bar{h}_t(\pi), \qquad \phi \in [0, \pi). \quad (37)$$

Using (35), we find that

$$\mu(A_t(0) \cap A_t(\phi))$$
$$= \begin{cases} -4 \sum_{k \text{ odd}} \frac{\gamma_k^t}{k} \sin\left(k\frac{\phi}{2}\right) + 2\pi \sum_{k \text{ odd}} \gamma_k^t - 2\pi \sum_{k \text{ even}} \gamma_k^t, & \text{if } \phi \in [0, \pi] \\ 4 \sum_{k \text{ odd}} \frac{\gamma_k^t}{k} \sin\left(k\frac{\phi}{2}\right) + 2\pi \sum_{k \text{ odd}} \gamma_k^t - 2\pi \sum_{k \text{ even}} \gamma_k^t, & \text{if } \phi \in [-\pi, 0]. \end{cases}$$

The result is now obtained by deriving a Fourier expansion of the latter expression and comparing with (36). $\square$

***Example 12.*** In the particular case where $g(s) = 1$ and

$$\bar{h}_t(\phi) = h_t(\phi) = \gamma_0^t + \gamma_1^t \cos \phi,$$

we find that

$$\lambda_0^t = \left[2\pi - \frac{16}{\pi}\right]\gamma_1^t - 2\pi\gamma_0^t,$$

$$\lambda_j^t = \frac{16}{\pi}\frac{1}{(2j)^2 - 1}\gamma_1^t, \qquad j = 1, 2, \ldots.$$

It follows that

$$(\lambda_j^t)^{-1} = \alpha_t + \beta_t j^2, \qquad j = 1, 2, \ldots, \tag{38}$$

where $\alpha_t = -\pi/(16\gamma_1^t)$ and $\beta_t = \pi/(4\gamma_1^t)$. Under the assumption of a normal Lévy basis, (38) is a special case of the $p$th order model considered in Jónsdóttir and Jensen (2005) with $p = 1$ and $\alpha$ proportional to $\beta$. Note that requirements (34) imply that $\gamma_0^t = t - \gamma_1^t$ and $\gamma_1^t > 0$. It does not seem to be possible to obtain $p$th order models with $p > 1$ using this approach.

## 6. An application to tumour growth

In Schmiegel (2006), snapshots of a growing brain tumour in vitro were analyzed using the approach described in this paper; see Figure 10. The data were first studied in Brú et al. (1998).

A detailed initial analysis of the covariance structure showed negative spatial covariances and a need for modelling both small and large scale fluctuations in the growth process. The model used was an exponential spatio-temporal Lévy model of the form

$$R_t(\phi) = \exp\biggl\{\mu_t + \alpha(t)\int_{t-T(t)}^{t-t_0(t)}\int_{-\pi}^{\pi}\cos(\phi-\theta)Z(\mathrm{d}s\,\mathrm{d}\theta) \\ + \beta(t)\int_{t-t_0(t)}^{t}\int_{\phi-h_t(s-t+t_0(t))}^{\phi+h_t(s-t+t_0(t))} Z(\mathrm{d}s\,\mathrm{d}\theta)\biggr\}. \tag{39}$$

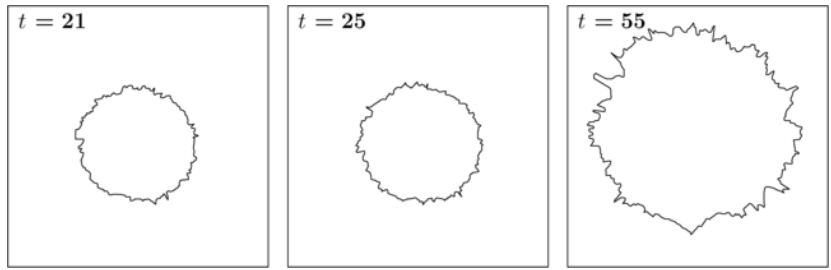

**Figure 10.** Profiles of a growing brain tumour in vitro at time points $t = 21$, 25 and 55.

**Table 2.** The estimated parameters for model (39), using a Gaussian Lévy basis

| $t$ | $T(t)$ | $t_0(t)$ | $\alpha(t)$ | $\beta(t)$ | $\phi_0(t)$ |
|---|---|---|---|---|---|
| 21 | 21 | 19 | 0.04 | $-0.033$ | 0.19 |
| 25 | 25 | 17 | 0.02 | $-0.033$ | 0.19 |
| 55 | 18 | 4 | 0.01 | $-0.067$ | 0.23 |

Here, $h_t$ is a deterministic and monotonically decreasing function defined on $[0, t_0(t)]$, satisfying $h_t(t_0(t)) = 0$ and $h_t(0) = \phi_0(t)/2$. Accordingly, the weight function is of the form

$$f_t(\xi;\phi) = \alpha(t)\cos(\phi-\theta)\mathbf{1}_{[t-T(t),t-t_0(t)]}(s) + \beta(t)\mathbf{1}_{[t-t_0(t),t]}(s)\mathbf{1}_{[0,h_t(s-t+t_0(t))]}(|\phi-\theta|).$$

The associated ambit set is shown in Figure 11. In Schmiegel (2006), a Gaussian Lévy basis was used and the function $h_t$ was assumed to be of the form

$$h_t(s) = \frac{\phi_0(t)}{2} - \frac{\phi_0(t)}{2t_0(t)}s, \qquad s \in [0, t_0(t)].$$

The parameters of model (39) were estimated by the method of moments, using the results given in Section 3. The estimated parameters are given in Table 2 and a simulation under the model with a Gaussian Lévy basis is shown in Figure 12.

Here, we will study the use of Gamma and inverse Gaussian Lévy bases. Simulations under the latter basis are shown in Figure 13. The inverse Gaussian Lévy basis is chosen such that $\mathbb{E}(R_t(\phi))$ and $\mathbb{V}(R_t(\phi))$ are the same as in the case where a Gaussian basis is used. The upper row of Figure 13 shows simulations where $\eta = 316$ and the lower row shows a simulation where $\eta = 5$. For $\eta = 316$, the inverse Gaussian Lévy basis provides fits of a similar quality as the normal basis, but for $\eta = 5$, more outbursts are observed, as is the case for the data. The difference is due to the fact that the inverse Gaussian distribution has heavier right tails for the latter choice of parameters. It should be noted that all of the profiles simulated under model (39) using the Lévy basis mentioned in this section show slightly more fluctuations on a local scale than the observed profiles. At present, we do not know whether this feature is caused by non-perfect model selection and estimation of parameters or artefacts due to the discretization in the simulation procedure.

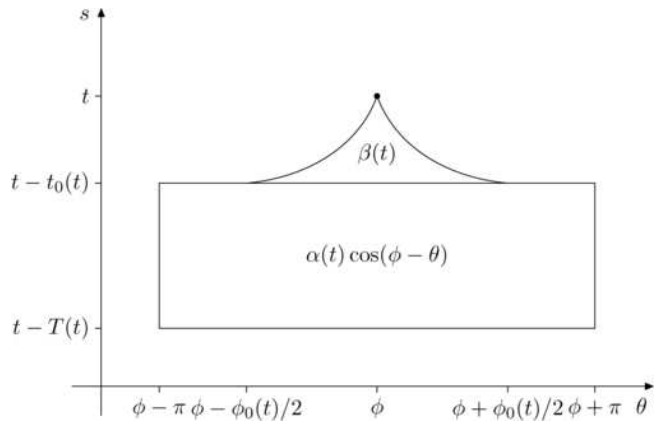

**Figure 11.** The ambit set $A_t(\phi)$ for the model defined by (39).

## 7. Discussion

In the present paper, we have given a condensed review of Lévy-based spatio-temporal modelling and shown its potential use in stochastic geometry and spatial statistics by developing Lévy-based growth models and space–time covariances on the circle. Below, we discuss further perspectives and topics for future research.

### 7.1. Related growth models

In the growth literature, there is a variety of growth models for objects in discrete space (cf., e.g., Bramson and Griffeath (1981), Qi *et al.* (1993), Lee and Cowan (1994), Kansal *et al.* (2000) and references therein). An important early example is the Richardson model, introduced in Richardson (1973). Here, the growth is described by a Markov process. For a growing object in the plane, the state at time $t$ is a random subset $Y_t$ of $\mathbb{Z}^2$ consisting of the 'infected sites'. An uninfected site is transferred to an infected site with a rate proportional to the number of infected nearest neighbours. It can be shown that if $Y_0$ consists of a single site, then $Y_t/t$ has a non-random shape as $t \to \infty$. Note that the growth model described in (20) of the present paper may be regarded as a continuous analogue of the Richardson model.

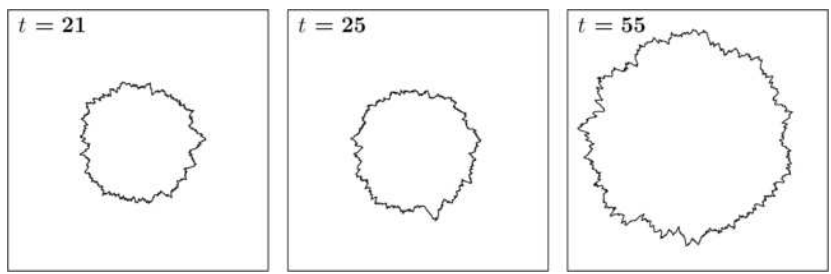

**Figure 12.** Simulation of model (39) for time points $t = 21$, 25 and 55, using a Gaussian Lévy basis.

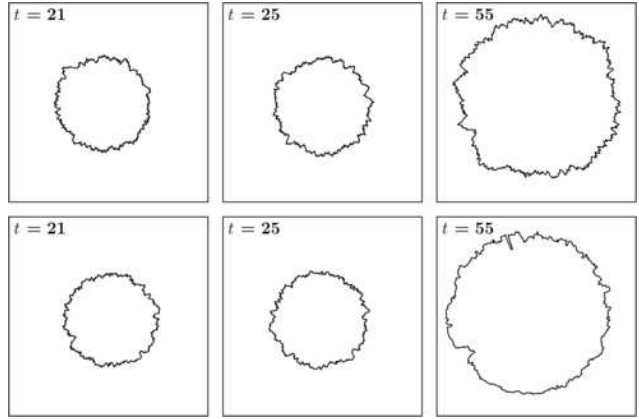

**Figure 13.** Simulations of model (39) for time points $t = 21$, 25 and 55 using an inverse Gaussian Lévy basis with $\eta = 316$ (upper row) and $\eta = 5$ (lower row).

A related growth model in continuous space has recently been discussed in Deijfen (2003). For planar objects, the model is constructed from a spatio-temporal Poisson point process on $\mathbb{R}^3$,

$$\Psi = \{(x_i, t_i)\}.$$

The random growing object $Y_t \subset \mathbb{R}^2$ is a subset of

$$\bigcup_{\{i \,:\, t_i \leq t\}} B(x_i, r),$$

constructed such that $Y_t$ is always connected. Here, $B(x, r)$ is a circular disc with centre $x$ and radius $r$. In this model, $t_i$ is thought of as a time point of outburst and $x_i$ as the location of the outburst in the tumour, say. A closely related discrete-time Markov growth model has been discussed in detail in Cressie and Hulting (1992). This model can be characterized as a sequence of Boolean models,

$$Y_{t+1} = \cup \{B(x_i, r) : x_i \in Y_t\},$$

where $\{x_i\}$ is a homogeneous Poisson point process in $\mathbb{R}^2$; see also Cressie and Laslett (1987) and Cressie (1991a, 1991b).

An issue of interest in growth modelling is the asymptotic shape of the growing object (cf., e.g., Durrett and Liggett (1981) and Deijfen (2003)). It is expected that it is also possible to obtain asymptotic results for Lévy-based growth models using the fact that Lévy bases are independently scattered random measures.

### 7.2. Lévy-driven Cox processes

Another interesting application of spatio-temporal Lévy models in spatial statistics is provided by Lévy-driven Cox processes (cf. Hellmund (2005), Hellmund *et al.* (2007) and Prokešová *et al.* (2006)). As an example, we may use $\exp(X_t(\sigma))$, where $X_t(\sigma)$ is given in (8), as driving random intensity for a spatio-temporal Cox process. If $Z$ is a Gaussian Lévy basis, then the resulting Cox process is log-Gaussian. Another example is obtained by assuming that $Z$ is a positive Lévy basis, the Cox process driven by $X_t(\sigma)$ then being a spatio-temporal shot-noise Cox process. Cox processes have been studied intensively in recent years (cf. Brix (1998), Wolpert and Ickstadt (1998), Brix (1999), Brix and Diggle (2001), Brix and Møller (2001), Brix and Chadoeuf (2002), Møller (2003)).

### 7.3. Estimation of model parameters

It still remains to develop inference procedures for Lévy-based growth models (and, more generally, for Lévy-based spatio-temporal models). There are several interesting problems here, including nonparametric estimation of the ambit sets.

A Fourier expansion of the radial function may be useful when making inference about the shape of the growing object (cf., e.g., Alt (1999) and Jónsdóttir and Jensen (2005)). Let us consider the Fourier coefficients of $R_t(\phi)$,

$$A_k^t = \frac{1}{\pi}\int_{-\pi}^{\pi} R_t(\phi)\cos(k\phi)\,\mathrm{d}\phi, \qquad B_k^t = \frac{1}{\pi}\int_{-\pi}^{\pi} R_t(\phi)\sin(k\phi)\,\mathrm{d}\phi,$$

$k = 0, 1, \dots$. Under the assumptions of Proposition 7, it can be shown that

$$A_k^t = \int_{-\pi}^{\pi}\int_{t-T(t)}^{t} a_k^t(s)\cos(k\theta)Z(\mathrm{d}\theta\,\mathrm{d}s), \qquad B_k^t = \int_{-\pi}^{\pi}\int_{t-T(t)}^{t} a_k^t(s)\sin(k\theta)Z(\mathrm{d}\theta\,\mathrm{d}s),$$

so the Fourier coefficients also follow a linear spatio-temporal Lévy model. It can be shown that, for $k \neq j$, $t, t' \geq 0$,

$$\mathrm{Cov}(A_k^t, A_j^{t'}) = \mathrm{Cov}(B_k^t, B_j^{t'}) = \mathrm{Cov}(A_k^t, B_j^{t'}) = 0$$

and

$$\mathrm{Cov}(A_k^t, A_k^{t'}) = \mathrm{Cov}(B_k^t, B_k^{t'}) = \tau_k(t, t'),$$

where $\tau_k(t, t')$ is given in Proposition 7.

In the case where $Z$ is a Gaussian Lévy basis, this means that $\{A_k^t\}_{t\in\mathbb{R}}$ and $\{B_k^t\}_{t\in\mathbb{R}}$, $k = 0, 1, \dots$, are independent Gaussian stochastic processes with covariance functions $\tau_k(t, t')$. If one observes $A_k^t$ and $B_k^t$ for some time points $t = t_1, \dots, t_n$ and some orders $k = 1, \dots, K_t$, the likelihood function is tractable.

# Acknowledgements

The authors thank Ole E. Barndorff-Nielsen for sharing his ideas with us which led to a more satisfactory solution to the original problem of analyzing the tumour growth data. This work was supported in part by grants from the Danish Natural Science Research Council and the Carlsberg Foundation.